

\documentstyle{amsppt}
\magnification 1200
\NoBlackBoxes
\def\ls{\vskip.25in}
\def\ss{\vskip.10in}
\def\C{\Bbb C}
\def\opl{\operatornamewithlimits{\oplus}}
\def\vep{\varepsilon}

\topmatter

\title{\bf Thickening Calabi-Yau moduli spaces}\endtitle
\ls
\author Z. Ran \endauthor
\affil Department of Mathematics\\ University of California\\
Riverside, CA  92521, USA\\ ziv\@ucrmath.ucr.edu \endaffil
\thanks  ${}^*$ Supported in part by NSF DMS 9202050
\endthanks
\abstract{\it We describe a kind of deformation of the anti-DeRham
algebra on a Calabi-Yau manifold $X$.  These are in 1-1 correspondence
with the total cohomology $\oplus H^i (X, \C)$. }
\endabstract
\endtopmatter
\document
\baselineskip=18pt

In his article [W] in this volume's precursor, Witten proposed, as a
possible approach to constructing a mirror map, a certain extended moduli
space $\Cal N$, a thickening of the usual moduli space $\Cal M$ of
complex structures on a Calabi-Yau manifold $X$.  Witten's proposal for
$\Cal N$ is couched in terms of Conformal Field Theory on $X$, and
it is not immediately obvious to the author to what extent it is or can
be made rigorous.  In any case, it appears to be an intriguing problem of
pure deformation theory to construct Witten's extended moduli
directly in terms of ordinary complex geometry, making no appeal to
Physics.
\ss
The purpose of this note is to describe a type of generalized or `exotic'
deformation of complex structure on a Calabi-Yau manifold $X$.
Our construction shares with Witten's the property that the set of first-order
deformations coincides with the total DeRham cohomology $\oplus H^i (X)$,
with the classical deformations embedded as $H^{1, n-1} (X)$; but the
author's incompetence in physical matters precludes any further comment
here as to any real connection between the two constructions.
\ss
Our approach focuses on the (Holomorphic) anti-DeRham or {\it Mahr Ed} algebra
$$
\Omega^{-\cdot}_X = \oplus \Lambda^i \Theta_X \ .
$$
This is a sheaf of graded anticommutative associative algebras.
Interestingly, when $X$ is endowed with a holomorphic volume form,
$\Omega^{-\cdot}_X$ not only becomes (additively) isomorphic with the
DeRham complex $\Omega^{\cdot}_X$ but also acquires a (graded)
bracket operation, known as the Schouten bracket,
extending the usual Lie bracket on $\Theta_X$.  This Schouten
bracket, which is familiar in symplectic geometry and mechanics,
turns $\Omega^{-\cdot}_X$ essentially into a graded Lie
algebra.  Moreover via the adjoint action, $\Omega^{-\cdot}$ can
essentially be identified as the differential graded Lie algebra
of derivations on the (associative) Mahr Ed algebra $(\Omega^{-\cdot},\wedge)$.
A suitable graded version of the usual connection between derivations
and deformations then yields the realization of the De Rham cohomology,
i.e., the hypercohomology of $(\Omega^{-\cdot}, d)$, as deformations of
the algebra $(\Omega^{-\cdot}, \wedge)$.
\ss
It must be stated that this work is not in final form with numerous
basic questions as yet unanswered and even unasked.  This applies
especially to the higher-order theory.  We hope to return to this
elsewhere.
\ls
\subheading{1.  The Schouten Lie algebra}
\ss
Fix an $n$-dimensional complex manifold $X$, not necessarily compact
(analogous considerations apply in the real case as well).  We denote as usual
by $\Theta_X$ and $\Omega_X$ its tangent and cotangent sheaves respectively,
and by $\Omega^{\cdot}_X = \opl\limits_{0}^n \Lambda^i \Omega_X$ the
DeRham algebra.  Our focus, however, will be on the latter's dual,
the {\it Mahr Ed} algebra
$$
\Omega^{-\cdot}_X = \opl\limits_{0}^n \Lambda^i \Theta_X .
$$
This is a sheaf of negatively graded, anticommutative associative
algebras.  Now it is an interesting observation going back to Schouten [S]
and Nijenhuis [N] (see also Lichnerowicz [L] and Koszul [K]), that the
Mahr Ed algebra carries an additional structure which is essentially
that of a {\it graded Lie algebra}.  To be precise, there is a bracket
operation [\quad] on $\Omega^{-\cdot}$, the {\it Schouten bracket},
with the following properties (where a,b,c denote homogenous elements):
\ss
{\parindent=20pt
\item{0.}  [\quad] coincides with the usual Lie bracket on $\Omega^{-\cdot}_X
= \Theta_X$;

\item{1.}  $\deg[a,b] = \deg \ a + \ \deg \ b \ + 1$;

\item{2.}  $[a,b] = - \vep(a,b)[b,a]$ where $\vep(a,b) = (-1)^{(\deg\ a +
1)(\deg\
b+1)}$;

\item{3.}  a Jacobi identity is satisfied:
$$
\vep(a,c) [a,[b,c]] + \vep(c,b)[c,[a,b]] + \vep (b,a)[b,[c,a]] = 0
$$

\item{4.}  $[a,\cdot]$ is a derivation of degree $\deg\ a+ 1$ on the
algebra $\Omega^{-\cdot}$.  (In fact, properties 0,1,2,4 already characterize
[\quad])
}

To exploit this structure, define the {\it Schouten Lie algebra} of $X$ by
$$
L^{\cdot}_X = \opl\limits_{0}^{\dim X-1} L^{-i}_X := \opl\limits_0^{\dim X-1}
\Lambda^{i+1} \Theta_X ,
$$
with the above Schouten bracket [\quad].  Then the above properties show that
$(L^{\cdot}_X, [\quad])$ forms a sheaf of graded Lie algebras on $X$.
Moreover $L^{\cdot}_X$ comes with a faithful graded Lie representation
$$
L^{\cdot}_X \to \text{Der}_{\C} (\Omega^{-\cdot}_X, \Omega^{-\cdot}_X ) .
\tag1.1
$$
Actually, this representation is surjective, i.e., an isomorphism,
though we won't need this fact.

\def\Io{I\llap{\raise1pt\hbox{$\circ$}}}
Now suppose given a volume form $\Phi$ on $X$, i.e., a nowhere vanishing
section of $\Omega^n_X$.  Then interior multiplication by $\Phi$ induces an
additive isomorphism
$$
i_{\Phi}: \Omega^{-\cdot}_X [n] \to \Omega^{\cdot}_X .
$$
via $i_{\Phi}$, the exterior derivative operator $d$ on $\Omega^{\cdot}_X$
may be imported over to $\Omega^{-\cdot}_X$, yielding an operator
$\delta =\delta_{\Phi}$ of degree $+1$ with $\delta ^2= 0$.  Note however
that as $\Io$ is not multiplicative, $\delta$ loses the all-important
derivation property of $d$.  In fact, the failure of $\delta$ to be a
derivation is measured precisely by the Schouten bracket:  this is because
of the following formula due to Koszul:
$$
\gather
[a,b]=(-1)^{\deg\ a} (\delta(ab)-\delta(a)b-(-1)^{\deg\ a \deg\ b}
\delta(b)a), \\
a,b\in \Omega^{-\cdot}\ \text{homogeneous}
\tag1.2
\endgather
$$
One consequence of Koszul's formula (1.2) is that the restriction of
$\delta$ on the Schouten algebra $L^{\cdot}$ may be identified with the
commutator with the operator $\delta$ on $\Omega^{-\cdot}$:
$$
\split
\delta(a) = [a,\delta] := \delta\circ a - (-1)^{\deg\ a+1} a \circ \delta\\
a \in L^{\cdot}\ \ \text{homogeneous}
\endsplit
\tag1.3
$$
In particular, we have
$$
\delta\big|_{L^{0}} \sim \ [\quad, d]\big|_{\Theta_X} = 0  .
\tag1.4
$$
(of course, this is what makes it sensible to restrict $\delta$ on $L^{\cdot}$
in the first place).  Note that the RHS of (1.1) a priori carries a
differential
that is commutator with $\delta$.  Thus the significance of (1.3) is precisely
that
the representation (1.1) is one of {\it differential} graded Lie algebras.
In particular, the cohomology of $(L^{\cdot},\delta)$ should be related to
properties of $\Omega^{-\cdot}$ as algebra, such as deformations.  This is
the basic idea of our approach.

To tie this in later with Dolbeault cohomology, note that the Schouten bracket
extends to a pairing
$$
A^{\circ,q} (L^{\cdot}) \times \Omega^{-\cdot} \to A^{0,q} (\Omega^{-\cdot})
=: A^{-\cdot,q}
$$
where $A^{\cdot,\cdot} (E)$ denotes real forms of type $(\cdot,\cdot)$ with
coefficients in $E$.  Naturally this is compatible with $\bar\partial$.

Finally it may be noted that via $i_{Phi}$, the pairing on
$\Omega^{\cdot}$ corresponding to exterior product on $\Omega^{-\cdot}$ is just
the so-called Yukawa pairing $*$.  Hence by (1.2) the Schouten bracket on
$\Omega^{\cdot}$ may be written as
$$
[\alpha,\beta] = \pm (d(\alpha * \beta)-d(\alpha)*\beta \pm d(\beta)*\alpha) .
$$
\bigskip

\subheading{2.  Exotic deformations}
\ss
Now fix a volume form $\Phi$ on our complex manifold $X$.   Motivated by
(1.4), we replace the Schouten algebra $L^{\cdot}_X$ by its subalgebra
$\hat{L}^{\cdot}_X$, the {\it restricted} Schouten algebra, defined by
$$
\align
\hat{L}^{0} &= \{ v \in L^{0} : L_v \Phi = 0\} , \quad L_v =
\text{\rm Lie derivative};\\
\hat{L}^{-i} = L^{-i} ,\ i > 0 .
\endalign
$$
Note the quasi-isomorphism
$$
\hat{L}^{\cdot} \sim \Omega^{\cdot}_X [n-1] .
$$
Combining this with the representation of $\hat{L}^{\cdot}$ as the DGL
algebra of volume-preserving derivations of the Mahr Ed algebra
$\Omega^{-\cdot}$ should, by general principles, give rise to a realization
of the De Rham cohomology $\Bbb H^{\cdot}(\Omega^{\cdot})$  as a suitable
kind of deformations of the `algebra with operator' $\Omega^{-\cdot}_X$,
which may be viewed as a sort of `exotic' or non-classical deformations
of the manifold $X$ itself.  We proceed to elaborate this idea.

We begin by defining a generalization of the usual notion of total complex
associated to a bicomplex.  Let $K^{\cdot,\cdot}$ be a doubly-indexed
array (not necessarily a bicomplex) of (abelian) objects with maps
$d', d''$ of bidegrees (1,0) and (0,1) respectively.  To this
we may associate as usual the {\it total object} $t^{\cdot}(K^{\cdot,\cdot})$
defined by
$$
t^i (K^{\cdot,\cdot}) = \opl\limits_{p+q=i} K^{p,q}
$$
By a {\it complexification} of $K^{\cdot,\cdot}$ we mean a structure  of
{\it complex} on $t^{\cdot} (K^{\cdot,\cdot})$, i.e. differentials $d: t^i
(K^{\cdot\cdot} \to t^{i+1} (K^{\cdot\cdot})$ with $d\circ d = 0$,
which are compatible with the original $d', d''$ in the sense that the
map $K^{p,q} \to K^{p`,q`}$, induced on subquotients by $d$, coincides
with $d'$ for $(p', q') = (p+1,q)$ and with $(-1)^q d''$ for
$(p', q') = (p, q+1)$, i.e., representing $d$ by a matrix, its
`near diagonals' should consist of $d'$ and $\pm d''$.  Of course,
if $K^{\cdot\cdot}$ happens to be a bicomplex, its usual associated
total complex is a complexification in the above sense, called the
{\it standard} one.  More generally, given a sub-bicomplex $J^{\cdot,\cdot}
\subseteq K^{\cdot,\cdot}$, a complexification is said to be
{\it standard on} $J^{\cdot,\cdot}$ if it preserves $t^{\cdot}
(J^{\cdot,\cdot}) \subseteq t^{\cdot}(K^{\cdot,\cdot})$ and agrees with the
standard on $t^{\cdot}(J^{\cdot,\cdot})$

\proclaim{Theorem}  There is a natural bijection between $\Bbb H^i
(F^{-j}\hat{L}^{-\cdot}) = \Bbb H^{i+n-1} (F^{n-1-j}\Omega^{\cdot}_X)$
and the set of sets of data as follows (up to a natural equivalence):

\item{(i)}  a surjection $E^{-\cdot} \to \Omega^{-\cdot}$ of sheaves of graded
anticommutative associative algebras, which is an isomorphism in degrees
$<j-n$, and which fits in a diagram with exact rows
$$
\matrix
&&&&&&&&&&E^{-n}&\simeq&\Omega^{-n}\\
&&&&&&&&&&\vdots&&\vdots\\
&&&&&&&&&&\downarrow&&\downarrow\\
0&\to&\Omega^{-n}&\to&A^{-n,0}&\to&\ldots&&A^{-n,i+j-2}
&\to&E^{j-n}&\to&\Omega^{j-n}&\to&0\\
&&\vdots&&&&&&&&\vdots&&\vdots\\&&\delta
\downarrow&&&&&&&&\downarrow&&\downarrow\delta\\
0&\to&\Omega^{-j}&\to&A^{-j,0}&\to&\ldots&\to&A^{-j, i+j-2}&\to&E^0&\to&\Cal
O&\to&0
\endmatrix
\tag2.1
$$

\item{(ii)}  a complexification $t^{\cdot} (A^{\cdot\cdot}, E^{\cdot})$ of the
middle portion $(A^{\cdot\cdot}, E^{\cdot})$ of {\rm (2.1)},
standard on $A^{\cdot\cdot}$ which fits in an exact sequence of complexes
$$
0\to\Omega^{-\cdot}/F^{-j+1}\Omega^{-\cdot} \to
t^{\cdot}(A^{\cdot\cdot},E^{\cdot})
\to \Omega^{-\cdot} \to 0 .
$$
\endproclaim

Let us give a brief indication how the above data are constructed out of
an element $u \in \Bbb H^i (F^{-j} \hat{L}^{-\cdot})$, thought of as
a mixed \v Cech-Dolbeault cohomology class.  First, $u$ gives
rise to a 1-cocycle $(a_{\alpha\beta})$ with values in ${}_{\bar{\partial}}
A^{-j,i+j-1} (L^-)$,  the $\bar\partial$-closed forms and from the
action of the latter on $\Omega^{-\cdot}$ we get the $E^{-\cdot}$.  Next
we have a 0-cochain $(b_\alpha )$ with values in $A^{-j+1, i+j-1} (L^{-\cdot})$
such that $b_\alpha - b_\beta = \partial \bar{a}_{\alpha\beta}$.  This
yields the differentials $\delta$ on $E^{\cdot}$ in (2.1).
However, we do not necessarily have a complex yet, i.e. $\delta^2\neq 0$.
The additional data of complexification is then obtained by pulling
$(b_\alpha)$ gradually further up the Hodge filtration in the
group $\Bbb H^i (F^{-j} \hat{L}^{-\cdot})$.

Note that aside from ordinary vector fields, the Schouten algebra
contains only elements of strictly negative degree, which are {\it nilpotent}
derivations of index $\leq n$.  On the other hand, corresponding results
in the classical case indicate for $X$ Calabi-Yau that all the
deformations involving $\hat{L}^0$ are unobstructed, i.e., extend to
all orders.  Thus our moduli space is an infinitesimal thickening of a
smooth space with tangent space essentially $\oplus H^i (\Theta_X) =
\oplus H^{n-1,i} (X)$.   In particular for $n=3$ the reduced moduli
has tangent space $H^{2,1} \oplus H^{2,2}$, with the two pieces supposedly
interchanged by mirror symmetry.

Note that, as a special case of the above, the deformations
corresponding to $F^{n-2} H^n(X) = \Bbb H^1 (F^{-1}\hat{L}^{\cdot\cdot})$
correspond to diagrams
$$
0 \to \Omega^{-\cdot}[-1]_{\leq 0} \to A^{-\cdot} [-1]_{\leq 0} \to E^{\cdot}
\to \Omega^{-\cdot} \to 0
$$
with $E^{\cdot}$ an algebra with operator $\delta$ satisfying
$\delta^2 = 0$; i.e. to extensions
$$
[D^{\cdot}] \in \ \text{Ext}_{\text{alg}}^1 (\Omega^{-\cdot}, {}_{\bar\partial}
A^{-\cdot,1} [-1]_{\leq 0} )
$$
where the Ext is to be interpreted in the above sense.
\bigskip
\subheading{3.  $H^1$ and local systems}
\ss
To help develop some feeling for the above exotic deformations, we
now consider in more detail what is probably the simplest non-classical
case, that of $H^1_{DR}$ (refusing to be bothered by the fact that for
most---and by some definition, all---Calabi-Yau manifolds this
group actually vanishes!).  We do not require $X$ to be compact or
K\"ahlerian.

Now in the case of $H^1_{DR} (X) = \Bbb H^{-n+2} (\hat{L}^{-\cdot}_X)$, the
construction of \S2 simplifies considerably.  It suffices to look at the
quotient algebra $E^0 \oplus E^{-1}$ of $E^{\cdot}$, together with the
differential $\delta : E^{-1} \to E^0$.  These fit in an exact diagram
$$
\matrix
0&\to&\Cal O&\to&E^{-1}&\to&\Omega^{-1}&\to&0\\
&&\downarrow&&\downarrow&&\downarrow\\
0&\to&\Omega&\to&E^0&\to&\Cal O&\to&0
\endmatrix
\tag3.1
$$
where we have identified $\Cal O \simeq \Omega^{-n}, \Omega \simeq
\Omega^{1-n}$.  Roughly, it is the $H^{0,1}$ part that yields
$E^0\oplus E^{-1}$, while the $H^{1,0}$ part yields the differential
$\delta$: in particular, when the $H^{0,1}$ part is zero, i.e.
$E^0\oplus E^{-1}$ is the trivial extension, the essential part of
$\delta$ is an operator $\Omega^{-1}\to\Omega$ of the form $u\mapsto
\omega\delta(u)$ for some $\omega\in H^0(\Omega)$.  Incidentally, a
similar simplification holds for any $F^{i-1} H^i_{DR} (X)$.

Now on the other hand, recall the usual interpretation of
$H^1_{DR}$ in terms of local systems:  there is a natural bijection
between $H^1_{DR} (X)$ and local systems, i.e. locally constant
sheaves $\Cal L$ on $X$, which fit in an exact sequence
$$
0\to\Bbb C_X \to \Cal L \to \Bbb C_X \to 0 .
\tag3.2
$$
Here the $H^{0,1}$ part of the data corresponds to any of the vector
bundles $F^i :=\Cal L \opl\limits_{\Bbb C_X} \Omega^i_X$, viewed as
extensions, while the $H^{1,0}$ part corresponds to an integrable
connection
$$
F^0 \to F^1 = F^0 \otimes \Omega^1_X,
$$
or, for that matter, to any of the operators $F^i \to F^{i+1}$ induced
by exterior derivative on $\Omega^{\cdot}_X$.

The question then is how to go back and forth directly between the
data (3.1) and (3.2).  The key to this comes from the diagram
$$
\matrix
0&\to&\Cal O_X&\to&E^{-1}&\to&\Omega^{-1}&\to&0\\
&&\|&&\downarrow&&\downarrow\\
0&\to&\Cal O_X&\to&F^{-1}&\to&\Cal O_X&\to&0\\
&& \downarrow&&\downarrow&&\| \\
0&\to&\Omega_X&\to&E^0&\to&\Cal O_X&\to&0\\
&&\| &&\downarrow&&\downarrow\\
0&\to&\Omega_X &\to&F^0&\to&\Omega_X&\to&0 .
\endmatrix
\tag3.3
$$
Now given $(E^{\cdot},\delta)$ we may consider the module of
differentials $\Omega_{E^\cdot}$, which inherits a grading from
$E^{\cdot}$, and put
$$
A^{\cdot} = (\Omega_{E^{\cdot}})^{\cdot} \otimes_{E^{\cdot}}
\Omega^{\geq -1} .
$$
We thus have a natural derivation $E^{\cdot} \overset
d_{E^{\cdot}}\to\rightarrow \Omega_{E^\cdot} \to A^{\cdot}$.  By universalit of
$d_{E^\cdot}$ the
operator $\delta : E^{-1} \to E^0$ factors through $A^{-1}$,
so we get maps
$$
E^{-1} \to A^{-1} \to E^0\to A^0 ,
$$
and clearly this sequence fits in a diagram like (3.3).  We may then simply
set
$$
\Cal L = \ker (A^{-1} \to A^0)
$$
and check easily that this is the appropriate local system.

To construct $(E^{\cdot}, \delta)$ from $\Cal L$ we use a variant of a
method used earlier, e.g. in [R].  Consider the sheaves $B^i = \Omega^i \oplus
\Cal L \otimes \Omega^{i+1}$ which we view as filtered:
$$
\Omega^{i+1} \subset \Cal L \otimes \Omega^{i+1} \subset B^i
\tag3.4
$$
For $j=0,-1$ define $E^j$ as the pairs $(\varphi, u)$ where $u\in \Omega^j$ and
$\varphi \in \text{Hom}_{\Cal O} (B^i, B^{i+j})$ is filtration-preserving,
induces multiplication by $u$ on $\Cal L \otimes \Omega^{i+1}$ and
multiplication by $u+\delta u$ on $\Omega^{i+1}\oplus \Omega^i =
B^i/\Omega^{i+1}$ (this is
independent of $i$).  The operator $\delta: E^{-1} \to E^0$ is defined
by
$$
\delta (\varphi, u) = ([\varphi, d], \delta u
)
$$
where $d: B^i \to B^{i+1}$ is the obvious differential.  It is not hard
to check (compare [R]) that $(E^{\cdot},\delta)$ is a deformation as
in (3.1) and is the correct one corresponding to $\Cal L$ as above.
\bigskip
\centerline{\bf References}
\ss
\item{[K]}  Koszul, J.L.:  ``Crochet de Schouten-Nijenhuis et cohomologie".
Ast\'erisque (hors s\'erie) (1985), 257--271.

\item{[L]}  Lichrerowicz, A.:  ``Las vari\'et\'es de Poisson et leurs
alg\`ebres de Lie associ\'es".  J. Differential Geometry
{\bf 12} (1977), 253-300.

\item{[N]}  Nijenhuis, A.:  ``Jacobi-type identities for bilinear
differential concomitants of certain tensor fields".  Indag. math
{\bf 17} (1955), 390-403.

\item{[R]}  Ran, Z.:  ``Canonical infinitesimal deformations".  Preprint

\item{[S]}  Schouten, J.A.:  ``On the differential operators of first
order in tensor calculus", in:  Convegno Internat. Geom. Diff.
(Italy 1953), pp. 1-7.  Roma:  Cremonese 1954.

\item{[W]}  Witten, E.:  ``Mirror manifolds and topological field
theory", in:  S.T. Yau, ed., ``Essays on mirror manifolds",
pp. 120-160.  Hong Kong:  International Press 1992.

\enddocument